\documentclass{amsart}
\usepackage{amssymb}
\usepackage{latexsym}
\usepackage{enumerate}

\newtheorem{thm}{Theorem}[section]
\newtheorem{defin}[thm]{Definition}
\newtheorem{cor}[thm]{Corollary}
\newtheorem{prop}[thm]{Proposition}
\newtheorem{lemma}[thm]{Lemma}
\newtheorem{eg}[thm]{Example}
\newtheorem{rmk}[thm]{Remark}

\newcommand{\naturals}{\ensuremath{\mathbb{N}}}
\newcommand{\integers}{\ensuremath{\mathbb{Z}}}
\newcommand{\rationals}{\ensuremath{\mathbb{Q}}}
\newcommand{\complex}{\ensuremath{\mathbb{C}}}

\newcommand{\pf}[1][]{\emph{Proof#1:\ }}

\newcommand{\abs}[1]{\ensuremath{\left|#1\right|}}
\newcommand{\defsetshort}[1]{\ensuremath{\left\{#1\right\}}}
\newcommand{\defsetspan}[1]{\ensuremath{<\!\!#1\!\!>}}
\newcommand{\bincoeff}[2]{\ensuremath{\left(\begin{array}{c}#1\\#2\end{array}\right)}}
\newcommand{\idmatrix}[1]{\ensuremath{\mbox{I}_{#1}}}
\newcommand{\idelement}[1]{\ensuremath{\mbox{Id}_{#1}}}

\renewcommand{\mod}[1]{\ensuremath{\ \left(\mbox{mod }#1\right)}}

\newcommand{\cohrank}[3]{\ensuremath{h^{#1}\left(#2,#3\right)}}
\newcommand{\cohgroup}[3]{\ensuremath{H^{#1}\left(#2,#3\right)}}

\newcommand{\symgroup}[1]{\ensuremath{\mathbb{S}_{#1}}}
\newcommand{\altgroup}[1]{\ensuremath{\mathbb{A}_{#1}}}
\newcommand{\dihgroup}[1]{\ensuremath{\mathbb{D}_{#1}}}
\newcommand{\cyclgroup}[1]{\ensuremath{\integers_{#1}}}
\newcommand{\glgroup}[2]{\ensuremath{\mbox{GL}_{#1}\left(#2\right)}}
\newcommand{\slgroup}[3][]{\ensuremath{\mbox{SL}^{#1}_{#2}\left(#3\right)}}
\newcommand{\pglgroup}[2]{\ensuremath{\mbox{PGL}_{#1}\left(#2\right)}}
\newcommand{\pslgroup}[2]{\ensuremath{\mbox{PSL}_{#1}\left(#2\right)}}

\newcommand{\autgroup}[2][]{\ensuremath{\mbox{Aut}_{#1}\left(#2\right)}}
\newcommand{\proj}[2][]{\ensuremath{\mathbb{P}_{#1}^{#2}}}

\newcommand{\deffun}[2]{\ensuremath{#1\rightarrow#2}}
\newcommand{\deffunname}[3]{\ensuremath{#1:#2\rightarrow#3}}

\newcommand{\funcKer}[1]{\ensuremath{\mbox{ker}\left(#1\right)}}
\newcommand{\funcIm}[1]{\ensuremath{\mbox{Im}\left(#1\right)}}
\newcommand{\funcOrb}[2][]{\ensuremath{\mbox{Orb}_{#1}\left(#2\right)}}

\title{Weakly-exceptional quotient singularities in prime dimension}
\author[D.\,Sakovics]{Dmitrijs Sakovics}
\address{Fakult{\"a}t f{\"u}r Mathematik\\Universit{\"a}t Wien\\Oskar-Morgenstern-Platz $1$\\$1090$ Wien\\Austria}
\email{dmitrijs.sakovics@univie.ac.at}

\begin{document}

\begin{abstract}
 A singularity is said to be weakly-exceptional if it has a unique purely log terminal blow up. This is a natural generalization of the surface singularities of types $D_{n}$, $E_{6}$, $E_{7}$ and $E_{8}$.  Since this idea was introduced, quotient singularities of this type have been classified in dimensions up to $5$. This paper looks at such singularities in dimension $p$, where $p$ is an arbitrary prime number.
\end{abstract}
\maketitle

\section{Introduction}

The central notion of this paper was first mentioned in V.~Shokurov's paper on flips (see~\cite{Shokurov92-Eng}), where (among other things) he looked at the properties of exceptional divisors that appear during (partial) blow-ups of the $2$-dimensiona $A$-$D$-$E$ singularities. This analysis has since been generalised as follows:
\begin{defin}
Let $\left(V\ni O\right)$ be a germ of a Kawamata log terminal singularity. The singularity is said to be exceptional if for every effective \rationals-divisor $D_V$ on the variety $V$, such that the log pair $\left(V,\ D_V\right)$ is log canonical, there exists at most one exceptional divisor over the point $O$ with discrepancy $−1$ with respect to the pair $\left(V,\ D_V\right)$.
\end{defin}
\begin{thm}[{\cite[Theorem~3.7]{Cheltsov-Shramov11a}}]
Let $\left(V\ni O\right)$ be a germ of a Kawamata log terminal singularity. Then there exists a birational morphism \deffunname{\pi}{W}{V}, such that the following hypotheses are satisfied:
\begin{itemize}
 \item the exceptional locus of $\pi$ consists of one irreducible divisor $E$ such that $O\in\pi\left(E\right)$,
 \item the log pair $\left(W,\ E\right)$ has purely log terminal singularities.
 \item the divisor $−E$ is a $\pi$-ample \rationals-Cartier divisor.
\end{itemize}
Such a morphism is called a \emph{plt blow-up} of the singularity.
\end{thm}
\begin{defin}
 We say that a singularity $\left(V\ni O\right)$ is \emph{weakly-exceptional} if it has a unique plt blow-up.
\end{defin}
This paper seeks to extend the following example:
\begin{eg}[see~{\cite[Section~5.2.3]{Shokurov92-Eng}}]\label{eg:dim2}
 Consider the case of $N=2$. The singularities in this case follow the well-known $A$-$D$-$E$ classification, all of them being quotient singularities of the type $\complex^2/G$ (for different groups $G\subset\slgroup{2}{\complex}$. The singularities of type $A_n$ correspond to (reducible) lifts of cyclic groups $\cyclgroup{n+1}\subset\autgroup{\proj{1}}$ to \slgroup{2}{\complex}, singularities of type $D_n$ correspond to lifts of dihedral groups $\dihgroup{n-2}\subset\autgroup{\proj{1}}$, and the singularities $E_6$, $E_7$, and $E_8$ --- to lifts of $\altgroup{4}$, $\symgroup{4}$ and~$\altgroup{5}$ respectively. Note that in this case the group action is irreducible exactly when the singularity is of type $D$ or $E$. 
 
 Rephrasing~\cite[Section~5.2.3]{Shokurov92-Eng}, a singularity is exceptional exactly when it is of type $E_6$, $E_7$ or~$E_8$; it is weakly-exceptional exactly if it is of type $D_n$ or $E_n$. Note that here a singularity is weakly-exceptional if and only if the group action that gives rise to it is irreducible (and it is exceptional if and only if the action is primitive).
\end{eg}
The last observation turns out to be partially true in general:
\begin{thm}{{\cite[Proposition~2.1]{Prokhorov00}}}
 If $G\subset\slgroup{N}{\complex}$ gives rise to an exceptional singularity, then the action of $G$ is primitive.
\end{thm}
\begin{prop}
 If $G\subset\slgroup{N}{\complex}$ gives rise to an weakly-exceptional singularity, then the action of $G$ is irreducible.
\end{prop}
\pf Similar to that of \cite[Proposition~2.1]{Prokhorov00}. \qed
\begin{rmk}
 However, the reverses of the last two statements are not true: for instance, the group $\altgroup{5}\subset\slgroup{3}{\complex}$ has an irreducible primitive action, but the corresponding singularity is neither exceptional nor weakly-exceptional (see~\cite[end of Section~3]{Sakovics12}).
\end{rmk}

Since the $A$-$D$-$E$ singularities are all quotient singularities, it makes sence to look at the case of quotient singularities in higher dimensions too.
\begin{rmk}
 Consider the singularity $\complex^N/G$ (where $G\subset\glgroup{N}{\complex}$). The definitions above mean that its exceptionality is dependent not on the group $G$ itself, but on its image under the natural projection to \pglgroup{N}{\complex}. Thus, using the Chevalley--Shephard--Todd theorem (see~\cite[Theorem~4.2.5]{Springer77}), it is possible to simplify the problem by assuming that $G\subset\slgroup{N}{\complex}$ (by assuming that $G$ contains no pseudoreflections and then choosing a convenient lift of the image of $G$ in \pslgroup{N}{\complex}).
\end{rmk}
A number of papers have been written about exceptional and weakly-exceptional quotient singularities in higher dimensions. The exceptional quotient singularities have all been classified in dimensions $3$, $4$, $5$, $6$ and $7$ in \cite{Cheltsov-Shramov11a}, \cite{Cheltsov-Shramov11b} and~\cite{Markushevich-Prokhorov}. It is has also been conjectured that there are no exceptional quotient singularities in higher dimensions.

The weakly-exceptional quotient singularities have also been classified in dimensions $3$, $4$ and $5$ (see~\cite{Sakovics12} and~\cite{Sakovics14}). This paper generalises this classification to all other prime dimensions (i.e.\ dimensions $q$, where $q$ is a prime integer) by the following result:
\begin{thm}[Main theorem]\label{thm:prime:main}
 Let $q$ be a positive prime integer. Then there are at most finitely many finite irreducible subgroups $\Gamma\subset\slgroup{q}{\complex}$, such that the singularity of $\complex^q/\Gamma$ is not weakly-exceptional.
\end{thm}
Note that the same result cannot hold in non-prime dimensions (for a counterexample, see~\cite[Theorem~1.15]{Sakovics12}). However, it is hoped that in non-prime dimensions all the counterexamples can be put into a small number of families --- see Remark~\ref{rmk:counterexample:nonprime}.

\section{Preliminaries}
In order to prove the main theorem, several previously known results need to be considered. When trying to show that a quotient singularity is weakly-exceptional, the following results become very useful:
\begin{prop}\label{prop:general:subgroup}
 Let $G\subset\glgroup{N}{\complex}$ be a finite subgroup, and $H\subset G$ a subgroup. If the singularity of $\complex^N/G$ is not weakly-exceptional, then neither is the singularity of $\complex^N/H$.
\end{prop}
\pf Immediate from the definition.\qed
\begin{prop}\label{prop:semiinvars}
 Let $G\subset\slgroup{N}{\complex}$ be a finite subgroup with a semi-invariant of degree $d<N$. Then the singularity of $\complex^N/G$ is not weakly-exceptional.
\end{prop}
\pf Immediate consequence of~\cite[Theorem~3.15]{Cheltsov-Shramov11a}.\qed

\begin{thm}[{\cite[Theorem~1.12]{Cheltsov-Shramov11c}}]\label{thm:general:WE-invariants}
 Let $G$ be a finite group in \glgroup{n+1}{\complex} that does not contain reflections. If $\complex^{n+1}/G$ is not weakly-exceptional, then there is a $\bar{G}$-invariant, irreducible, normal, Fano type projectively normal subvariety $V\subset\proj{n}$ such that 
\[
 \deg{V}\leq\bincoeff{n}{\dim{V}}
\]
and for every $i\geq1$ and for every $m\geq0$ one has
\[
 \cohrank{i}{\proj{n}}{\mathcal{O}_{\proj{n}}\left(\dim{V}+1\right)\otimes\mathcal{I}_V} = \cohrank{i}{V}{\mathcal{O}_V\left(m\right)}=0,
\]
\[
 \cohrank{0}{\proj{n}}{\mathcal{O}_{\proj{n}}\left(\dim{V}+1\right)\otimes\mathcal{I}_V}\geq\bincoeff{n}{\dim{V}+1},
\]
where $\mathcal{I}_V$ is the ideal sheaf of the subvariety $V\subset\proj{n}$. Let $\Pi$ be a general linear subspace of \proj{n} of codimension $k\leq\dim{V}$. Put $X=V\cap\Pi$. Then $\cohrank{i}{\Pi}{\mathcal{O}_\Pi\left(m\right)\otimes\mathcal{I}_X}=0$ for every $i\geq0$ and $m\geq k$, where $\mathcal{I}_X$ is the ideal sheaf of the subvariety $X\subset\Pi$. Moreover, if $k=1$ and $\dim{V}\geq2$, then $X$ is irreducible, projectively normal, and $\cohrank{i}{X}{\mathcal{O}_X\left(m\right)}=0$ for every $i\geq1$ and $m\geq1$.
\end{thm}

Since thesee considerations rely on the group action used, one needs to define some terms common in the study of representation theory of finite groups:
\begin{defin}
 Given a group $G\subset\glgroup{N}{\complex}$, a \emph{system of imprimitivity} for $G$ is a set \defsetshort{V_1,\ldots,V_k} of subspaces of $\complex^N$, such that $\dim{V_i}>0$ $\forall i$, $V_i\cap V_j=\defsetshort{0}$ whenever $i\neq j$, $V_1\otimes\ldots\otimes V_k=\complex^N$, and for any $g\in G$ and $i\in\defsetshort{1,\ldots,k}$, there exists $j\left(g,i\right)\in\defsetshort{1,\ldots,k}$, such that $g\left(V_i\right)=V_{j\left(g,i\right)}$.
\end{defin}
\begin{rmk}
  Clearly, any group $G\subset\glgroup{N}{\complex}$ has at least one system of imprimitivity, namely, \defsetshort{\complex^N}.
\end{rmk}
\begin{defin}\label{def:prelim:prim}
 A group $G\subset\glgroup{N}{\complex}$ is \emph{primitive} if it has exactly one system of imprimitivity.
\end{defin}
This leads to the following well-known result:
\begin{lemma}[{Jordan's theorem --- see, for example, \cite{Frobenius11}}]\label{lemma:Jordan}
 For any given $N$, there are only finitely many finite primitive subgroups of \slgroup{N}{\complex}.
\end{lemma}
\begin{defin}\label{def:prelim:irred}
 A group $G\subset\glgroup{N}{\complex}$ is \emph{irreducible} if for any system of imprimitivity \defsetshort{V_1,\ldots,V_k} for $G$, the action of $G$ permutes the subspaces $V_1,\ldots,V_k$ transitively.
\end{defin}
\begin{prop}\label{prop:sysOimpr:irred-samedim}
 If a group $G\subset\glgroup{N}{\complex}$ with a system of imprimitivity \defsetshort{V_1,\ldots,V_k} is irreducible, then $k$ divides $N$, and \[\dim{V_1}=\dim{V_2}=\cdots=\dim{V_k}=N/k.\]
\end{prop}
\pf Since $G$ is irreducible, given $i,j\in\defsetshort{1,\ldots,k}$, there exists $g_{i,j}\in G$ such that $g_{i,j}\left(V_i\right)=V_j$. Therefore, $\dim{V_i}=\dim{V_j}$. Applying this for different pairs $(i,j)$, get $\dim{V_1}=\dim{V_2}=\cdots=\dim{V_k}=d$, some $d\in\integers$. Since the pairwise intersections between the $V_i$-s are trivial, and they span all of $\complex^N$, $kd=N$.\qed
\begin{defin}\label{def:monomial}
 A group $G\subset\glgroup{N}{\complex}$ is \emph{monomial} if there exists a system of imprimitivity \defsetshort{V_1,\ldots,V_k} for $G$, such that for all $i\in\defsetshort{1,\ldots,k}$, $\dim{V_i}=1$.
\end{defin}
These groups have additional structure that can be exploited by using the following results:
\begin{prop}\label{prop:prime:mono-prim}
 Let $q>1$ be a prime number, and let $G\subset\slgroup{q}{\complex}$ be a finite irreducible subgroup. Then the action of $G$ is either primitive or monomial.
\end{prop}
\pf Given any system of imprimitivity for $G$, Proposition~\ref{prop:sysOimpr:irred-samedim} implies that all the subspaces in that system must have the same dimension $d$, with $d|q$. Since $q$ is prime, $d\in\defsetshort{1,q}$. If there exists a system with $1$-dimensional subspaces, then the action of $G$ is monomial. Otherwise, the action of $G$ must be primitive.\qed

\begin{prop}\label{prop:monomialStructure}
  Let $G\subset\glgroup{N}{\complex}$ be a finite monomial subgroup. Then have $G\cong D\rtimes T$, where $D$ is abelian, and $T\subseteq\symgroup{N}$. Given a system of imprimitivity \defsetshort{V_1,\ldots,V_N} for this group and choosing $0\neq x_i\in V_i$ for every $i$, the set $\defsetshort{x_{1},\ldots,x_{N}}$ forms a basis for $\complex^N$. In this basis, every element of $D$ is a diagonal matrix, and for every element $g\in G\setminus D$, there exists some $i,j\in\defsetshort{1,\ldots,N}$ with $i\neq j$ and $g\left(x_i\right)\in V_j$.
\end{prop}
\pf
 Since $G$ is monomial, it has at least one system of imprimitivity \defsetshort{V_1,\ldots,V_k}, such that all the $V_i$-s have dimension $1$. Since $V_1,\ldots,V_k$ span $\complex^N$, must have $k=N$. The action of $G$ permutes $V_1,\ldots,V_N$, so have a homomorphism \deffunname{\pi}{G}{\symgroup{N}} defined by these permutations. Let $D=\funcKer{\pi}\unlhd G$ and $T=\funcIm{\pi}\subseteq\symgroup{N}$. Clearly, $G=D\rtimes T$.

 For every $i$, pick a non-zero element $x_i\in V_i$. Since $V_i$ is one-dimensional, $x_i$ spans $V_i$, and so $\defsetshort{x_{1},\ldots,x_{N}}$ is a basis for $\complex^N$. Given any $d\in D$, $d\left(V_i\right)=V_i$ for every $i$, and so $d$ must be a diagonal matrix in the chosen basis. Therefore, $D$ is abelian.

 Let $g\in G$, such that $g\left(x_i\right)\in V_i$ for all $i$. Then $\pi\left(g\right)$ is the trivial permutation in \symgroup{N}, and so $g\in\funcKer{\pi}=D$. So for any $g\in G\setminus D$, there exist~$i\neq j$ with $g\left(x_i\right)\in V_j$.
\qed
\begin{prop}\label{prop:monomial:irred-trans}
 Let $G\subset\glgroup{N}{\complex}$ be a finite monomial subgroup, and let $G\cong D\rtimes T$ be the decomposition from Proposition~\ref{prop:monomialStructure}. If $G$ is irreducible, then $T$ is transitive.
\end{prop}
\pf
 Assume $T\subseteq\symgroup{N}$ is not transitive. Let $x_1$ be a basis vector from Proposition~\ref{prop:monomialStructure}. Consider the subspace $V$ of $\complex^N$ spanned by \funcOrb[G]{x_1}. Since $T$ is not transitive, there exists $j\in\defsetshort{1,\ldots,N}$ such that $j\not\in\funcOrb[T]{1}$. Therefore, $V_j\cap V=\defsetshort{0}$, and so $V\neq\complex^N$. However, by construction $V$ must be $G$-invariant, and so $G$ is not irreducible.
\qed

\begin{lemma}[{\cite[\S8.1]{Serre77}}]\label{lemma:monomialReps}
 If $A$ is an abelian normal subgroup of a group $G$, then the degree of each irreducible representation of $G$ divides the index $\left(G:A\right)$ of $A$ in $G$.
\end{lemma}

Finally, several miscellaneous results will be needed in the technical part of the proof:
\begin{thm}[\cite{Lam-Leung00}]\label{thm:misc:sums_of_roots}
 Given $m\in\integers$, $m>1$, let $W\left(m\right)$ be the set of integers $n\geq0$, for which there exist $\omega_1,\ldots,\omega_n\in\complex$ with $\omega_i^m=1\ \forall i$ and $\omega_1+\cdots+\omega_n=0$. Take $m=p_1^{a_1}\cdots p_r^{a_r}$ the prime decomposition of $m$. Then
\[
 W\left(m\right)=\naturals p_1+\naturals p_2+\cdots+\naturals p_r
\]
\end{thm}

\begin{defin}[\cite{Davis79}]
 An $n$-by-$n$ matrix $M$ is called \emph{circulant} if it is of the form
\[
M=\left(\begin{array}{ccccc}
          a_{1}&a_{2}&\cdots&a_{n-1}&a_{n}\\
          a_{n}&a_{1}&\cdots&a_{n-2}&a_{n-1}\\
          \vdots&\vdots&\ddots&\vdots&\vdots\\
          a_{3}&a_{4}&\cdots&a_{1}&a_{2}\\
          a_{2}&a_{3}&\cdots&a_{n}&a_{1}\\
         \end{array}
\right)
\]
 for some numbers $a_1,\ldots,a_n\in\complex$.
\end{defin}
\begin{lemma}[{\cite[\S3.2]{Davis79}}]\label{lemma:misc:circulant_eigenvalues}
 For any circulant matrix $M$ with $a_1,\ldots,a_n$ as above and any $\omega\in\complex$ with $\omega^n=1$, $M$ has an eigenvector $v=\left(1,\omega,\omega^2,\ldots,\omega^{n-1}\right)^T$ with eigenvalue $\lambda=\sum_{i=1}^{n}{a_i\omega^{i-1}}$. All the eigenvalues of $M$ are of this form.
\end{lemma}
\pf It is easy to check that vectors of this form are indeed eigenvectors of $M$ with relevant eigenvalues. These form a set of $n$ linearly independent eigenvectors (can be seen via Theorem~\ref{thm:misc:sums_of_roots}), so these are all the possible eigenvalues and eigenvectors for $M$.\qed

\begin{lemma}\label{lemma:prime:det-nonzero}
 Let $q\in\naturals$ be prime, and consider the following $q\times q$ matrix with integer coefficients:
\[
M=\left(\begin{array}{ccccc}
          a_{1}&a_{2}&\cdots&a_{q-1}&a_{q}\\
          a_{q}&a_{1}&\cdots&a_{q-2}&a_{q-1}\\
          \vdots&\vdots&\ddots&\vdots&\vdots\\
          a_{3}&a_{4}&\cdots&a_{1}&a_{2}\\
          a_{2}&a_{3}&\cdots&a_{q}&a_{1}\\
         \end{array}
\right)
\]
Assume that $a_i\geq0$ $\forall i$, and $0<d=\sum_{i=1}^{q}{a_i}<q$. Then the determinant of $M$ is not zero.
\end{lemma}
\pf Consider the matrix $M$ over \complex, and assume $\det{M}=0$. Then one of the eigenvalues of $M$ must be zero. So, by Lemma~\ref{lemma:misc:circulant_eigenvalues},
\[
 a_1+\omega a_2+\omega^2a_3+\ldots+\omega^{q-1}a_q=0
\]
for some $\omega$ with $\omega^q=1$. Since all the $a_i$-s are non-negative integers, this is a sum of exactly $d=\sum_{i=1}^{q}{a_i}$ $q$-th roots of unity. So  $d\in W\left(q\right)$. However, by Theorem~\ref{thm:misc:sums_of_roots} and using the fact that $q$ is prime, $W\left(q\right)=\naturals q$. But by the initial assumption, $0<d<q$, producing a contradiction.\qed

\section{Proof of main result}

The aim of this section is to prove Theorem~\ref{thm:prime:main}. From now on, assume that $q\geq3$ is a prime, and $\Gamma\subset\slgroup{q}{\complex}$ be a finite irreducible subgroup, such that the singularity of $\complex^q/\Gamma$ is not weakly-exceptional.

By Jordan's theorem (see Lemma~\ref{lemma:Jordan}) there are only finitely many primitive finite subgroups of \slgroup{q}{\complex}. Therefore, for the purposes of this proof, one can assume that the group $\Gamma$ is imprimitive. Furthermore, $q$ is assumed to be prime, so (by Proposition~\ref{prop:prime:mono-prim}) $\Gamma$ must be monomial.

\begin{lemma}\label{lemma:prime:getTau}
 Assume $G\subset\slgroup{q}{\complex}$ is a finite irreducible monomial subgroup. Setting $G\cong D\rtimes T$ as in Proposition~\ref{prop:monomialStructure}, there exists $\tau\in G\setminus D$ and a basis $e_{1},\ldots,e_{q}$ for $\complex^q$, such that $\tau^q=\idelement{G}$, and $\tau$ acts by
\[
 \tau\left(e_i\right)=e_{i+1}\ \forall i<q;\ \tau\left(e_q\right)=e_1
\]
\end{lemma}
\pf 
 Since $G$ is irreducible, $T$ must be a transitive subgroup of \symgroup{q} (by Proposition~\ref{prop:monomial:irred-trans}) and must thus contain a cycle of length $q$ (since $q$ is prime). Take $\tau\in\Gamma$, such that $\pi\left(\tau\right)$ is a generator of this cycle. Let $e_1\in V_1$ be a non-zero vector. Then, renaming the $V_i$-s if necessary, $\tau^i\left(e_1\right)\in V_{i+1}$ (for $1\leq i<q$). Set $e_i=\tau^{i-1}\left(e_1\right)$ ($2\leq i\leq q$). Clearly, $\tau\left(e_{q}\right)=\alpha e_1$ for some $\alpha\in\complex$.

 Since all the subspaces $V_i$ are disjoint and one-dimensional, $e_i$ must generate $V_i$, and so $e_1,\ldots,e_q$ must form a basis for $\complex^q$. Also, since $g\in D=\ker{\pi}$, and $\tau$ permutes the subspaces $V_i$ non-trivially, $\tau\not\in D$. Since $\tau\in G\subseteq\slgroup{q}{\complex}$ and $q$ odd, one also observes that $\alpha=1$, and so $\tau$ acts as stated above.\qed

\begin{cor}\label{cor:prime:getSmallerGroup}
 There exists a subgroup $G=D\rtimes\cyclgroup{q}\subseteq\Gamma$ generated by $D$ and $\tau$. The singularity of $\complex^q/G$ is not weakly-exceptional, and $\abs{\Gamma}\leq \left(q-1\right)!\abs{G}$.
\end{cor}
\pf
 Take $G$ generated by $D$ and the element $\tau\in\Gamma$ obtained in Lemma~\ref{lemma:prime:getTau}. Clearly, $G\subseteq\Gamma$ and, looking at the action of $\tau$, $G\cong D\rtimes\cyclgroup{q}$. Since $G\subseteq\Gamma$, the singularity of $\complex^q/G$ is not weakly-exceptional by Proposition~\ref{prop:general:subgroup}. Finally, $\Gamma\subseteq D\rtimes\symgroup{q}$ (by Proposition~{\ref{prop:monomialStructure}), so
\[\abs{\Gamma}\leq\frac{\abs{\symgroup{q}}}{\abs{\cyclgroup{q}}}\abs{G}=\left(q-1\right)!\abs{G}\]
\qed

From now on, fix the group $G$ constructed above, the abelian normal subgroup $D\lhd G$, the element $\tau\in G$ and the basis $e_{1},\ldots,e_{q}$ for $\complex^q$ constructed in Lemma~\ref{lemma:prime:getTau}.

It is now advantageous to obtain a specialised criterion for determining whether or not a group of this form gives rise to a weakly-exceptional singularity.

\begin{prop}\label{prop:prime:RepDimensions}
 Any irreducible representation of $G$ (given above) over $\complex$ is either $1$-dimensional or $q$-dimensional.
\end{prop}
\pf Lemma~\ref{lemma:monomialReps} implies that $\left(G:D\right)=q$, which is a prime.\qed

\begin{lemma}[generalising~{\cite[Theorem~3.4]{Cheltsov-Shramov11c}}]\label{lemma:prime:semiinvar}
 Let $q$ be an odd prime and assume $G\subset\slgroup{q}{\complex}$ is a finite imprimitive subgroup isomorphic to $A\rtimes\cyclgroup{q}$ for some abelian $A$. Then the singularity of $\complex^q/G$ is not weakly-exceptional if and only if $G$ has a (non-constant) semi-invariant of degree $d<q$.
\end{lemma}
\pf
 If $G$ does have a semi-invariant of degree at most $q-1$, then the singularity is not weakly-exceptional by Proposition~\ref{prop:semiinvars}. Suppose that $G$ does not have any such semi-invariants, but the singularity is not weakly-exceptional.

 Then, by Theorem~\ref{thm:general:WE-invariants}, there exists a $\bar{G}$-invariant irreducible normal Fano type variety $V\subset\proj{q-1}$, such that $\deg{V}\leq\bincoeff{q-1}{\dim{V}}$ and $\cohrank{i}{V}{\mathcal{O}_V\left(m\right)}=0$ $\forall i\geq1$ $\forall m\geq0$ 
 (where $\mathcal{O}_V\left(m\right)=\mathcal{O}_V\otimes\mathcal{O}_{\proj{q-1}}\left(m\right)$).

 Let $n=\dim{V}$. Then, since $G$ has no semi-invariants of degree less than $q$, have $n\leq q-2$. Let $\mathcal{I}_V$ be the ideal sheaf of $V$. Then
 \[
  \cohrank{0}{V}{\mathcal{O}_V\left(m\right)} = \cohrank{0}{\proj{q-1}}{\mathcal{O}_{\proj{q-1}}\left(m\right)} - \cohrank{0}{\proj{q-1}}{\mathcal{I}_V\left(m\right)}
 \]
 For instance, $\cohrank{0}{V}{\mathcal{O}_V}=1$.

 Take any $m\in\integers$ with $0<m<q$. Let $W_m=\cohgroup{0}{\proj{q-1}}{\mathcal{I}_V\left(m\right)}$. This is a linear representation of $G$, so $q|\dim{W_m}$ (by Proposition~\ref{prop:prime:RepDimensions}, as $G$ has no semi-invariants of degree $m<q$). Since $q|\cohrank{0}{\proj{q-1}}{\mathcal{O}_{\proj{q-1}}\left(m\right)}$,
 \[
  \cohrank{0}{V}{\mathcal{O}_V\left(m\right)}\equiv0\mod{q}
 \]

 Since $\cohrank{0}{V}{\mathcal{O}_V\left(t\right)}=\chi\left(V,\mathcal{O}_V\left(t\right)\right)$ for any integer $t\geq0$, there exist integers $a_0,\ldots,a_n$, such that
 \[
  \cohrank{0}{V}{\mathcal{O}_V\left(t\right)}=P\left(t\right)=a_nt^n+a_{n-1}t^{n-1}+\cdots+a_1t+a_0
 \]

 Consider $P\left(t\right)$ as a polynomial over \cyclgroup{q}. Since
 \[
  P\left(m\right)=\cohrank{0}{V}{\mathcal{O}_V\left(m\right)}\equiv0\mod{q}
 \]
 whenever $0<m<q$, $P\left(t\right)$ has at least $q-1$ roots over \cyclgroup{q}. But $\deg{P}\leq n\leq q-2$, so $P\left(t\right)$ must be the zero polynomial over \cyclgroup{q}. In particular, $a_0\equiv0\mod{q}$. On the other hand, $a_0=P\left(0\right)=\cohrank{0}{V}{\mathcal{O}_V}=1\not\equiv0\mod{q}$, leading to a contradiction.\qed

Now let $f\left(x_1,\ldots,x_q\right)$ be a semi-invariant of $G$ of degree $d<q$ as given in the above lemma. Using the chosen basis, let
\[
 m\left(x_1,\ldots,x_q\right)=x_1^{a_1}x_2^{a_2}\cdots x_q^{a_q}
\]
be a monomial contained in $f$ (for some $a_i\in\integers_{\geq0}$, not all zero). This means that $f$ must contain all the monomials in the $\tau$-orbit of $m$. Furthermore, $\sum_i{a_i}=d<q$ and $\sum_{i=0}^q{\lambda^i\tau^i\left(m\right)}$ is a semi-invariant of $G$ whenever $\lambda^q=1$. So, without loss of generality, can assume
\[
 f\left(x_1,\ldots,x_q\right)=\left[m+\lambda\tau\left(m\right)+\cdots+\lambda^{q-1}\tau^{q-1}\left(m\right)\right]\left(x_1,\ldots,x_q\right)
\]

This semi-invariant can now be exploited to obtain a bound for the possible size of $D$.

First, consider the possible cyclic subgroups of $D$. Lemma~\ref{lemma:prime:det-nonzero} makes it possible to make the following deductions:
\begin{lemma}\label{lemma:prime:boundCycleSize}
 Let $g\in D$, and let $n$ be the smallest positive integer, such that $g^n$ is a scalar matrix. Then $n<q^{2q+1}$.
\end{lemma}
\pf Assume $n>1$. Since $g\in G\subset\slgroup{q}{\complex}$, $g^n=\zeta_q\mbox{I}_q$, where $\zeta_q$ is a $q$-th root of $1$ and \idmatrix{q} is the identity matrix. Then, since all the elements of $D$ are diagonal matrices,
\[
g=\zeta_q^{\beta_0}\left(\begin{array}{ccc}
          \zeta_{n}^{\beta_1}&&\\
          &\ddots&\\
          &&\zeta_{n}^{\beta_q}
         \end{array}
\right)
\]
where $\beta_i\in\integers$, not all zero, with $0\leq\beta_i<n$ $\forall i>0$; $0\leq\beta_0<q$. Since $n$ was taken to be minimal, the highest common factor of \defsetshort{n,\beta_1,\ldots,\beta_q} is $1$.

Now consider the polynomial $f$ of degree $d<q$ described above. Since we know $g\in G$, $g\left(f\right)=\lambda f$ for some $\lambda\in\complex$. Since $g^{nq}=\idmatrix{q}$ and all the monomials are $g$-semi-invariant, $\lambda=\zeta_q^{\beta_0}\zeta_{n}^C$, some $C\in\integers$. This is equivalent to:
\begin{eqnarray*}
 C&\equiv&a_1\beta_1+a_2\beta_2+\cdots+a_{q-1}\beta_{q-1}+a_q\beta_q\ \mod{n}\\
 &\equiv&a_1\beta_2+a_2\beta_3+\cdots+a_{q-1}\beta_q+a_q\beta_1\ \mod{n}\\
 &\equiv&a_1\beta_3+a_2\beta_4+\cdots+a_{q-1}\beta_1+a_q\beta_2\ \mod{n}\\
 &&\ldots\\
 &\equiv&a_1\beta_q+a_2\beta_1+\cdots+a_{q-1}\beta_{q-2}+a_q\beta_{q-1}\ \mod{n}
\end{eqnarray*}
This can be rewritten as
\[
 M\left(\beta_1,\ldots,\beta_q\right)^T\equiv C\left(1,\ldots,1\right)^T\mod{n}
\]
where $M$ is the matrix from Lemma~\ref{lemma:prime:det-nonzero}). However, since $\sum_{i=1}^{q}{a_i}=d$, $M$ also satisfies
\[
 M\left(1,\ldots,1\right)^T=d\left(1,\ldots,1\right)^T
\]

Take $v=d\left(\beta_1,\ldots,\beta_q\right)^T-C\left(1,\ldots,1\right)^T$. By linearity, $Mv\equiv0\mod{n}$. Multiplying both sides by the adjugate matrix of $M$, get:
\begin{eqnarray*}
 \left(d\beta_1-C\right)\det{M}&\equiv&0\mod{n}\\
 \left(d\beta_2-C\right)\det{M}&\equiv&0\mod{n}\\
 \ldots\\
 \left(d\beta_q-C\right)\det{M}&\equiv&0\mod{n}
\end{eqnarray*}
Therefore,
\[
 d\beta_1\det{M}\equiv d\beta_2\det{M}\equiv\cdots\equiv d\beta_q\det{M}\equiv C\det{M}\mod{n}
\]
This implies that $g^{d\det{M}}$ is a scalar matrix. By assumption, $0<d<q$ (in \integers), and, by Lemma~\ref{lemma:prime:det-nonzero}, $\det{M}\neq0$ (over \integers), so $\abs{d\det{M}}=Kn$ for some positive integer $K$. Thus, $n\leq\abs{d\det{M}}\leq q\abs{det{M}}$.

Now look at the entries $M_{i,j}$ of the matrix $M$. Since $0\leq a_k\leq d<q$ for all $k$, $\abs{M_{i,j}}\leq d<q$. Thus,
\[
 n\leq q\abs{\det{M}}\leq q\left(q\max_{i,j}{\abs{M_{i,j}}}\right)^q<q^{2q+1}
\]
\qed
\begin{cor}\label{cor:prime:MaxCycleSize}
 Let $\cyclgroup{m}\subseteq D$. Then $m\leq q^{2q+2}$. 
\end{cor}
\pf Take $g$ a generator of $\cyclgroup{m}\subseteq D$. Then for some $n\leq q^{2q+1}$, $g^n$ is a scalar matrix in \slgroup{q}{\complex}. Therefore, $g^{qn}=\idelement{G}$. So 
\[
 m\leq qn\leq q\cdot q^{2q+1}=q^{2q+2}
\]
\qed
\begin{lemma}\label{lemma:prime:MaxCycleCount}
 Let $\left(\cyclgroup{m}\right)^k\subseteq D\subset G\subseteq\Gamma\subset\slgroup{q}{\complex}$. Then $k\leq q$.
\end{lemma}
\pf Let $g_1,\ldots,g_k$ be a minimal set of generators of $\left(\cyclgroup{m}\right)^k\subseteq D$. Then for every $i>1$, $g_i\not\in\defsetspan{g_1,\ldots,g_{i-1}}$. Let $\zeta_m$ be a primitive $m$-th root of $1$. Then all the $g_i$ are diagonal matrices with some powers of $\zeta_m$ as diagonal entries. But any matrix in \slgroup{q}{\complex} has exactly $q$ diagonal entries, so at most $q$ such $g_i$-s can be chosen. Therefore, $k\leq q$.\qed
\begin{cor}\label{cor:prime:LimitD}
 $D\subseteq \bigotimes_{i=0}^{q^{2q+2}}{\left(\cyclgroup{i}\right)^q}$.
\end{cor}
\pf Immediate from Corollary~\ref{cor:prime:MaxCycleSize} and Lemma~\ref{lemma:prime:MaxCycleCount}.\qed
\begin{thm}\label{thm:prime:finManyMono}
 Given $q>3$, there are at most finitely many finite irreducible monomial groups $\Gamma\subseteq\slgroup{q}{\complex}$, such that the singularity of $\complex^q/\Gamma$ is not weakly-exceptional.
\end{thm}
\pf Let $\Gamma$ be such a group. Then by Corollary~\ref{cor:prime:getSmallerGroup}, there exists $G\subseteq\Gamma$, such that $G\cong D\rtimes\cyclgroup{q}$ and $\abs{\Gamma}\leq\left(q-1\right)!\abs{G}$. By Corollary~\ref{cor:prime:LimitD}, $D\subseteq \bigotimes_{i=0}^{q^{2q+2}}{\left(\cyclgroup{i}\right)^q}$, so there are at most finitely many such group $D$. It follows that there are at most finitely many such groups $G$, and hence at most finitely many such groups $\Gamma$.\qed
\begin{rmk}
 The bounds used here for the possible sizes of the groups $D$, $G$ and $\Gamma$ are by no means effective. However, improving them would make the proofs in this section a lot more technically complicated, and would not provide much insight into the structure of these groups or significantly improve the main result of this paper.
\end{rmk}

This result provides the last step needed for the proof of the main theorem of this paper:

\pf[ of Theorem~\ref{thm:prime:main}] Let $q$ be any positive prime integer. If $q=2$, no such groups $\Gamma$ exists (by Example~\ref{eg:dim2}), so assume $q\geq3$. Then, by Proposition~\ref{prop:prime:mono-prim}, $\Gamma$ is either monomial or primitive. If $\Gamma$ is primitive, then it must be among a finite list of groups by Lemma~\ref{lemma:Jordan}. If $\Gamma$ is monomial, it must belong to a finite set of groups by Theorem~\ref{thm:prime:finManyMono}. Therefore, there are only finitely many such subgroups $\Gamma\subset\slgroup{q}{\complex}$. \qed

Note that the proof of this theorem also provides the means of enumerating all the imprimitive groups that $\Gamma$ can be isomorphic (or conjugate) to for any given $q$. This would rely on making a list of all the possible matrices $M$ form Lemma~\ref{lemma:prime:boundCycleSize} and computing (the prime factorisations of) their determinants.

Although the bounds seen in the proofs are very large, choosing a specific value of $q$ quickly produces a fairly short list for the possible isomorphism classes of $D$:
\begin{eg}
 If $q=7$, then $D\subseteq\cyclgroup{7}\times\left(\cyclgroup{n\cdot d}\right)^6$, where the values of~$n$ and~$d$ are as follows:
\begin{center}
\begin{tabular}{ll}
 $d=2$&$n=2^6$\\
 $d=3$&$n$ is one of $2^3$, $3^6$, $43$\\
 $d=4$&$n$ is one of $2^{12}$, $29$, $71$, $547$\\
 $d=5$&$n$ is one of $2^6$, $2^3\cdot71$, $5^6$, $13^2$, $29\cdot113$, $43$, $197$, $421$, $463$\\
 $d=6$&$n$ is one of $2^9$, $2^6\cdot3^6$, $2^3\cdot29$, $2^6\cdot43$, $13^2$, $29^2$, $29\cdot449$, $41^2$, $43\cdot71$,\\& $113$, $197$, $211$, $379$, $463$, $757$, $2689$.
\end{tabular}
\end{center}
\end{eg}
\pf By easy direct computation. \qed

A further computation (to obtain the possible actions of the group $T$ on $D\subset D\rtimes T=\Gamma$) can reduce this list even further. An example of such a computation can be seen for the case of $q=5$ in~\cite{Sakovics12}.

\begin{rmk}\label{rmk:counterexample:nonprime}
 It is easy to see that there is no hope to have the same result when the dimension $N$ is not a prime number. The easiest way to see counterexamples for arbitrarilly high dimension is to take $N=n^2$ and write the coordinates of $\complex^N$ as entries of an $n\times n$ matrix. This gives a map \deffunname{\iota}{\slgroup{n}{\complex}\times\slgroup{n}{\complex}}{\slgroup{N}{\complex}}, where the copies of \slgroup{n}{\complex} act by left and (transposed) right multiplication. It is easy to choose pairs of finite subgroups $A,B\subset\slgroup{n}{\complex}$, such that $\Gamma=\iota\left(A,B\right)$ acts irreducibly, and the action of $\Gamma$ clearly has a degree $n$ invariant (given by the matrix determinant). Clearly, one can choose infinitely many different suitable pairs $\left(A,B\right)$, yielding infinitely many such groups $\Gamma\subset\slgroup{N}{\complex}$.
 
 It is, however, hoped that, for any given dimension $N$, groups in the image of \deffun{\slgroup{a}{\complex}\times\slgroup{b}{\complex}}{\slgroup{ab}{\complex}} account for all the infinite families of groups that give rise to singularities that are not weakly-exceptional.
\end{rmk}


\begin{thebibliography}{10}

\bibitem{Cheltsov-Shramov11a}
I.~Cheltsov and C.~Shramov.
\newblock {On exceptional quotient singularities}.
\newblock {\em Geom. Topol.}, 15(4):1843--1882, 2011.

\bibitem{Cheltsov-Shramov11b}
I.~Cheltsov and C.~Shramov.
\newblock {Six-dimensional exceptional quotient singularities}.
\newblock {\em Math. Res. Lett.}, 18(6):1121--1139, 2011.

\bibitem{Cheltsov-Shramov11c}
I.~Cheltsov and C.~Shramov.
\newblock {Weakly-exceptional singularities in higher dimensions}.
\newblock {\em J. Reine Angew. Math.}, 689:201--241, 2014.

\bibitem{Davis79}
P.J. Davis.
\newblock {\em {Circulant matrices}}.
\newblock John Wiley \& Sons, New York-Chichester-Brisbane, 1979.
\newblock A Wiley-Interscience Publication, Pure and Applied Mathematics.

\bibitem{Frobenius11}
F.G. Frobenius.
\newblock {\em {{\"U}ber den von {L}. {B}ieberbach gefundenen Beweis eines
  Satzes von {C}. {J}ordan}}.
\newblock K{\"o}nigliche Akademie der Wissenschaften, 1911.

\bibitem{Lam-Leung00}
T.~Y. Lam and K.H. Leung.
\newblock {On vanishing sums of roots of unity}.
\newblock {\em J. Algebra}, 224(1):91--109, 2000.

\bibitem{Markushevich-Prokhorov}
D.~Markushevich and Yu.G. Prokhorov.
\newblock {Exceptional quotient singularities}.
\newblock {\em Amer. J. Math.}, 121(6):1179--1189, 1999.

\bibitem{Prokhorov00}
Yu.G. Prokhorov.
\newblock {Sparseness of Exceptional Quotient Singularities}.
\newblock {\em Mathematical Notes}, 68(5-6):664--667, 2000.

\bibitem{Sakovics12}
D.~Sakovics.
\newblock {Weakly-exceptional quotient singularities}.
\newblock {\em Cent. Eur. J. Math.}, 10(3):885--902, 2012.

\bibitem{Sakovics14}
D.~Sakovics.
\newblock {Five-dimensional weakly exceptional quotient singularities}.
\newblock {\em Proc. Edinb. Math. Soc.}, 57(1):269--279, 2014.

\bibitem{Serre77}
Jean-Pierre Serre.
\newblock {\em {Linear representations of finite groups}}.
\newblock Springer-Verlag, New York, 1977.
\newblock Translated from the second French edition by Leonard L. Scott,
  Graduate Texts in Mathematics, Vol. 42.

\bibitem{Shokurov92-Eng}
V.V. Shokurov.
\newblock {3-fold log flips}.
\newblock {\em Russ. Ac. SC Izv. Math.}, 40(1):95--202, 1993.

\bibitem{Springer77}
T.A. Springer.
\newblock {\em {Invariant theory}}.
\newblock {Lecture Notes in Mathematics, Vol. 585}. Springer-Verlag, Berlin-New
  York, 1977.

\end{thebibliography}
\end{document}